\documentclass{article}
\usepackage{epsfig}               % Please include any local package you are using
\usepackage{mathrsfs}
\usepackage{amssymb} %Extra symbols 
\usepackage{color}
\usepackage{hyperref}
\usepackage{url}
\definecolor{refcolor}{RGB}{0,0,190}
\hypersetup{
    colorlinks,
    citecolor=refcolor,
    filecolor=refcolor,
    linkcolor=refcolor,
    urlcolor=refcolor
}

\author{Ovidiu-Cristinel Stoica}               % Here insert the authors name and the address
\title{Spacetimes with Singularities
\footnote{``Simion Stoilow'' - Institute of Mathematics of the Romanian Academy,
  21 Calea Grivitei Street, 010702 Bucharest, Romania.
	Author's email: h\,o\,l\,o\,t\,r\,o\,n\,i\,x\,@\,g\,m\,a\,i\,l\,.\,c\,o\,m}
\footnote{Partially supported by Romanian Government grant PN II Idei 1187.}
\footnote{This article contains the talk named ``Spacetimes with Singularities'', at ``The 10th International Workshop on Differential Geometry and its Applications'' held between August 25 and August 30 in Constan\c{t}a, Rom{\^a}nia. It appeared in \cite{Sto12e}.}
}               % Here insert the title
\date{}               % this must be left empty

\begin{document}
\maketitle

\begin{abstract}
We report on some advances made in the problem of singularities in general relativity.

First is introduced the singular semi-Riemannian geometry for metrics which can change their signature (in particular be degenerate). The standard operations like covariant contraction, covariant derivative, and constructions like the Riemann curvature are usually prohibited by the fact that the metric is not invertible. The things become even worse at the points where the signature changes. We show that we can still do many of these operations, in a different framework which we propose. This allows the writing of an equivalent form of Einstein's equation, which works for degenerate metric too.

Once we make the singularities manageable from mathematical viewpoint, we can extend analytically the black hole solutions and then choose from the maximal extensions globally hyperbolic regions. Then we find space-like foliations for these regions, with the implication that the initial data can be preserved in reasonable situations. We propose qualitative models of non-primordial and/or evaporating black holes.

We supplement the material with a brief note reporting on progress made since this talk was given, which shows that we can analytically extend the Schwarzschild and Reissner-Nordstr\"om metrics at and beyond the singularities, and the singularities can be made degenerate and handled with the mathematical apparatus we developed.

\end{abstract}

% Here comes the body of the text
%\tableofcontents

%--------------------------------------------------------
% Defines

\newtheorem{theorem}{Theorem}[section]
\newtheorem{lemma}[theorem]{Lemma}
\newtheorem{proposition}[theorem]{Proposition}
\newtheorem{corollary}[theorem]{Corollary}
\newtheorem{definition}[theorem]{Definition}
\newtheorem{example}[theorem]{Example}
\newtheorem{remark}[theorem]{Remark}
   
\def\({\left(}
\def\){\right)}

\newcommand{\eqref}[1]{(\ref{#1})}

\newcommand{\R}{\mathbb{R}}
\newcommand{\N}{\mathbb{N}}
\newcommand{\C}{\mathbb{C}}

\newcommand{\ie}{\textit{i.e.} }
\newcommand{\cf}{\textit{cf.} }
\newcommand{\eg}{\textit{e.g.} }

\newcommand{\mf}[1]{\mathfrak{#1}}
\newcommand{\mc}[1]{\mathcal{#1}}
\newcommand{\ms}[1]{\mathscr{#1}}

\newcommand{\de}{\textnormal{d}}
\newcommand{\supp}[1]{\textnormal{supp}(#1)}

\newcommand{\metric}[1]{\langle#1\rangle}
\newcommand{\annihg}{\coannih{g}}
\newcommand{\idxannih}[2]{#1{}^{#2}{}}
\newcommand{\idxcoannih}[2]{#1{}_{#2}{}}
\newcommand{\radix}[1]{\idxcoannih{#1}{\circ}}
\newcommand{\annih}[1]{\idxannih{#1}{\bullet}}
\newcommand{\coannih}[1]{\idxcoannih{#1}{\bullet}}
\newcommand{\coradix}[1]{\idxannih{#1}{\circ}}

\newcommand{\vectmodule}{\mf X}
\newcommand{\fivect}[1]{\vectmodule(#1)}
\newcommand{\fiscal}[1]{\ms F(#1)}
\newcommand{\annihforms}[1]{\annih{\mc A}(#1)}
\newcommand{\discformsk}[2]{A_d{}^{#1}(#2)}
\newcommand{\fiformk}[2]{\mc A^{#1}(#2)}
\newcommand{\annihprod}[1]{\coannih{\langle\!\langle#1\rangle\!\rangle}}
\newcommand{\srformsk}[2]{\annih{\ms A}{}^{#1}(#2)}
\newcommand{\tensors}[3]{\mc T{}^{#1}_{#2}#3}

\newcommand{\cocontr}{{{}_\bullet}}
\newcommand{\lie}{\mc L}
\newcommand{\kosz}{\mc K}
\newcommand{\der}{\nabla}
\newcommand{\dera}[1]{\der_{#1}}
\newcommand{\derb}[2]{\dera{#1}{#2}}
\newcommand{\derc}[3]{({\derb{#1}{#2}})(#3)}
\newcommand{\lder}{\der^{\flat}}
\newcommand{\ldera}[1]{\lder_{#1}}
\newcommand{\lderb}[2]{\ldera{#1}{#2}}
\newcommand{\lderc}[3]{(\lderb{#1}{#2})(#3)}
\newcommand{\ric}{\textnormal{Ric}}

\newcommand{\image}[3]{
\begin{figure}
\begin{center}
\epsfig{file=#1,width=#2 in}
\end{center}
\caption{\small{\label{#1}#3}}
\end{figure}}

%--------------------------------------------------------
\section{Introduction}

%--------------------------------------------------------
\subsection{The universe as a geometry}

As geometers, we love to live, in our minds, in a geometric world. But how much of the physical world is geometry?

Many mathematicians, including Riemann, Hamilton, Clifford, tried to describe the physical reality as a geometric structure, and even pondered whether matter could actually be a geometric property of the space.

The first materialization of their idea is Einstein's \textit{general relativity}, which establishes the connection between the matter fields and geometry, by Einstein's equation
\begin{equation}
\label{eq_einstein_idx}
	G_{ab} + \Lambda g_{ab} = \frac{8\pi \mathcal G}{c^4} T_{ab}.
\end{equation}

The mathematical framework of general relativity is \textit{semi-Riemannian} (or pseudo-Riemannian) geometry, which is a generalization of Riemannian geometry to \textit{indefinite metrics}.

The predictions of general relativity were confirmed, and all the tests devised by physicists were passed successfully, proving that general relativity provides an accurate description of the world.

%--------------------------------------------------------
\subsection{The trouble with general relativity}

But two big clouds seem to be a menace for this paradise:
\begin{enumerate}
	\item 
An \textit{external} one -- the (apparent) incompatibility with \textit{quantum theory}, which describes other very important properties of the matter
	\item 
An \textit{internal} one -- the occurrence of singularities.
\end{enumerate}

This presentation is about a research program concerning the problem of singularities, and the progress made so far.
This problem comes from physics, but it is purely geometric.

%--------------------------------------------------------
\subsection{The singularity theorems}

A collapsing star can become a black hole - an object which bends the lightcones so that nothing, not even light, can escape. Inside the black hole there is a singularity.

When this was initially observed from the theoretical description of a collapsing star, it was hoped that this problem cannot actually occur, and it is an artifact of the perfect spherical symmetry.
But the singularity theorems show that, under very general assumptions, general relativity implies that such singularities occur with necessity  (Penrose and Hawking \cite{Pen65,Haw66i,Haw66ii,Haw67iii,HP70,HE95}). Subsequently, it was shown that the conditions leading to singularities are more general by Christodoulou \cite{Chr09}, and then even more general, by Klainerman and Rodnianski \cite{KR09}.

It was stated (\cite{HP70,Haw76,ASH91,HP96,Ash08,Ash09}) that

\vspace{3mm}
\textit{General relativity predicts its own breakdown!}
\vspace{3mm}

%--------------------------------------------------------
\subsection{Problems with the degenerate metric}
  
Let's take the simplest case of singularity in semi-Riemannian geometry: the metric becomes degenerate.
  
Then, the inverse of the metric, $g^{ab}$, which is so necessary for raising indices, and contracting between covariant indices, cannot be constructed.
  
The Christoffel symbol involves in its definition the inverse metric $g^{ab}$. The Riemann curvature needs $g^{ab}$, because its definition involves contractions.
Therefore, we cannot construct the Levi-Civita connection from the Koszul formula.

Actually, for the case where the signature of the metric is constant, Demir Kupeli constructed a \textit{Koszul derivative}, which unfortunately is not connection or covariant derivative, is not unique and is not invariant \cite{Kup87a,Kup87b,Kup87c,Kup96}. He used a (screen) subbundle of the tangent bundle, which is maximal so that the restriction of the metric to this bundle is non-degenerate.
  
He then constructed a curvature $R_\nabla$, depending on the connection. The quantity $\metric{R_\nabla(X,Y)Z,T}$ turned out to be a tensor field. To work free from the dependence of $\nabla$ and $R_\nabla$ on the screen bundle, Kupeli developed a part of his results in the quotient bundle $TM/\radix{T}M$ (we denoted by $\radix{T}M$ the subbundle of $TM$ made at each point $p\in M$ of the degenerate vectors from $T_pM$).

Unfortunately, his result won't apply for our case, because we need the metric to change its signature.

As we shall see, by a different approach than Kupeli's, we can construct some canonically (uniquely) defined invariants like the covariant derivative for an important class of covariant tensors (and differential forms), and the Riemann curvature tensor. This approach will allow us to work also with variable signature. In the case when the signature is constant, our Riemann curvature tensor turns out to be just Kupeli's $\metric{R_\nabla(X,Y)Z,T}$.

%--------------------------------------------------------
\section{Singular semi-Riemannian geometry}

We now proceed to an introduction in singular semi-Riemannian geometry, and our contribution which is intended to help solving the problems of singularities in general relativity. Except for the first section, the following material was developed in \cite{Sto11a,Sto11b,Sto11d}, where detailed proofs were given.

%--------------------------------------------------------
\subsection{Singular semi-Riemannian geometry}

\begin{definition}(see \eg \cite{Kup87b}, \cite{Pam03} p. 265 for comparison)
\label{def_sing_semiRiemm_man}
A \textit{singular semi-Riemannian manifold} is a pair $(M,g)$, where $M$ is a differentiable manifold, and $g\in \Gamma(T^*M \odot_M T^*M)$ is a symmetric bilinear form on $M$, named \textit{metric tensor} or \textit{metric}. If the signature of $g$ is fixed, then 
$(M,g)$ is said to be with \textit{constant signature}. If the signature of $g$ is allowed to vary from point to point, $(M,g)$ is said to be with \textit{variable signature}. If $g$ is non-degenerate, then $(M,g)$ is named \textit{semi-Riemannian manifold}. If $g$ is positive definite, $(M,g)$ is named \textit{Riemannian manifold}.
\end{definition}

\begin{definition}(\cf \eg \cite{Bej95} p. 1, \cite{Kup96} p. 3 and \cite{ONe83} p. 53)
Let $(V,g)$ be a finite dimensional inner product space, where the inner product $g$ may be degenerate. The totally degenerate space $\radix{V}:=V^\perp$ is named the \textit{radical} of $V$. An inner product $g$ on a vector space $V$ is non-degenerate if and only if $\radix{V}=\{0\}$.
\end{definition}

\begin{definition}(see \eg \cite{Kup87b} p. 261, \cite{Pam03} p. 263)
We denote by $\radix{T}M$ and we call \textit{the radical of $TM$} the following subset of the tangent bundle: $\radix{T}M=\cup_{p\in M}\radix{(T_pM)}$. We can define vector fields on $M$ valued in $\radix{T}M$, by taking those vector fields $W\in\fivect M$ for which $W_p\in\radix{(T_pM)}$.
\end{definition}

%--------------------------------------------------------
\subsection{Covariant contraction}

Here we make some preparations, and then introduce the covariant contraction, which is normally prohibited by the fact that the metric is degenerate, hence its inverse $g^{ab}$ -- which normally performs the covariant contraction -- cannot be defined \cite{Sto11a}.

\begin{definition}
We denote by $\annih{T}M$ the subset of the cotangent bundle defined as
\begin{equation}
	\annih{T}M=\bigcup_{p\in M}\annih{(T_pM)}
\end{equation}
where $\annih{(T_pM)} \subseteq T^*_pM$ is the space of covectors at $p$ which can be expressed as $\omega_p(X_p)=\metric{Y_p,X_p}$ for some $Y_p\in T_p M$ and any $X_p\in T_p M$. $\annih{T}M$ is a vector bundle if and only if the signature of the metric is constant. We can define sections of $\annih{T}M$ in the general case, by
\begin{equation}
	\annihforms{M}:=\{\omega\in\fiformk 1{M}|\omega_p\in\annih{(T_pM)}\textnormal{ for any }p\in M\}.
\end{equation}
\end{definition}

\begin{definition}
\label{def_co_inner_product}
On $\annih{T}M$ we can define a unique non-degenerate inner product $\annihg$ by $\annihg(\omega,\tau):=\metric{X,Y}$, where $X,Y\in\fivect M$, $\annih X=\omega$ and $\annih Y=\tau$. We alternatively use the notation $\annihprod{\omega,\tau}=\annihg(\omega,\tau)$.
\end{definition}

\begin{definition}
\label{def_radix_annih_tensor_field}
Let $T$ be a tensor of type $(r,s)$. We call it \textit{radical-annihilator} in the $l$-th covariant slot if  $T\in \tensors r{l-1}{M}\otimes_M\annih{T}M\otimes_M \tensors 0{s-l}{M}$.
\end{definition}

\begin{definition}
\label{def_contraction_covariant}
We can define uniquely the \textit{covariant contraction} or \textit{covariant trace} operator by extending the inner product $\annihg$. First, we define it on tensors $T\in\annih{T}M\otimes_M\annih{T}M$, by $C_{12}T=\annihg^{ab}T_{ab}$. This definition is independent on the basis, because $\annihg\in\annih{T}^*M\otimes_M\annih{T}^*M$. This operation extends by linearity to any tensors which are radical in two covariant indices.
We denote the contraction $C_{kl} T$ of a tensor field $T$ by 
\begin{equation}
T(\omega_1,\ldots,\omega_r,v_1,\ldots,\cocontr,\ldots,\cocontr,\ldots,v_s).
\end{equation}
\end{definition}

%--------------------------------------------------------
\subsection{Covariant derivative}

In the following we construct the replacement for the covariant derivative, for the case when the metric is allowed to become degenerate. It is not possible to just extend the definitions of the Levi-Civita connection and the covariant derivative, because they need the inverse of the metric, which is not defined. We will see that we can still do a lot without introducing them. Let's recall first the definition of the Koszul form and its properties, without proof, from \cite{Sto11a}.

\begin{definition}[The Koszul form]
\label{def_Koszul_form}
\textit{The Koszul form} is defined as
\begin{equation}
	\kosz:\fivect M^3\to\R,
\end{equation}
\begin{equation}
\label{eq_Koszul_form}
\begin{array}{llll}
	\kosz(X,Y,Z) &:=&\displaystyle{\frac 1 2} \{ X \metric{Y,Z} + Y \metric{Z,X} - Z \metric{X,Y} \\
	&&\ - \metric{X,[Y,Z]} + \metric{Y, [Z,X]} + \metric{Z, [X,Y]}\}.
\end{array}
\end{equation}
\end{definition}

\begin{theorem}
\label{thm_Koszul_form_props}
The Koszul form of a singular semi-Riemannian manifold $(M,g)$ has the following properties:
\begin{enumerate}
	\item \label{thm_Koszul_form_props_linear}
	It is additive and $\R$-linear in each of its arguments.
	\item \label{thm_Koszul_form_props_flinearX}
	It is $\fiscal M$-linear in the first argument:

	$\kosz(fX,Y,Z) = f\kosz(X,Y,Z).$
	\item \label{thm_Koszul_form_props_flinearY}
	Satisfies the \textit{Leibniz rule}:

	$\kosz(X,fY,Z) = f\kosz(X,Y,Z) + X(f) \metric{Y,Z}.$
	\item \label{thm_Koszul_form_props_flinearZ}
	It is $\fiscal M$-linear in the third argument:

	$\kosz(X,Y,fZ) = f\kosz(X,Y,Z).$
	\item \label{thm_Koszul_form_props_commutYZ}
	It is \textit{metric}:

	$\kosz(X,Y,Z) + \kosz(X,Z,Y) = X \metric{Y,Z}$.
	\item \label{thm_Koszul_form_props_commutXY}
	It is \textit{symmetric} or \textit{torsionless}:

	$\kosz(X,Y,Z) - \kosz(Y,X,Z) = \metric{[X,Y],Z}$.
	\item \label{thm_Koszul_form_props_commutZX}
	Relation with the Lie derivative of $g$:

	$\kosz(X,Y,Z) + \kosz(Z,Y,X) = (\lie_Y g)(Z,X)$.
	\item \label{thm_Koszul_form_props_commutX2Y}

	$\kosz(X,Y,Z) + \kosz(Y,Z,X) = Y\metric{Z,X} + \metric{[X,Y],Z}$.
	\end{enumerate}
for any $X,Y,Z\in\fivect M$ and $f\in\fiscal M$.
\end{theorem}

\begin{definition}[The lower covariant derivative]
\label{def_l_cov_der}
The \textit{lower covariant derivative} of a vector field $Y$ in the direction of a vector field $X$ is the differential $1$-form $\lderb XY \in \fiformk 1{M}$ defined as
\begin{equation}
\label{eq_l_cov_der_vect}
\lderc XYZ := \kosz(X,Y,Z)
\end{equation}
for any $Z\in\fivect{M}$.
The \textit{lower covariant derivative operator} is the operator
\begin{equation}
	\lder:\fivect{M} \times \fivect{M} \to \fiformk 1{M}
\end{equation}
which associates to each $X,Y\in\fivect{M}$ the differential $1$-form $\ldera XY$.
\end{definition}

\begin{definition}[see \cite{Kup96} Definition 3.1.3]
\label{def_radical_stationary_manifold}
A singular semi-Riemannian manifold $(M,g)$ is \textit{radical-stationary} if it satisfies the condition 
\begin{equation}
\label{eq_radical_stationary_manifold}
		\kosz(X,Y,\_)\in\annihforms M,
\end{equation}
for any $X,Y\in\fivect{M}$.
\end{definition}

\begin{definition}
\label{def_cov_der_covect}
Let $(M,g)$ be a radical-stationary semi-Riemannian manifold. We define the covariant derivative of a radical-annihilator $1$-form $\omega\in\annihforms{M}$ in the direction of a vector field $X\in\fivect{M}$ by
\begin{equation}
	\der:\fivect{M} \times \annihforms{M} \to \discformsk 1 M
\end{equation}
\begin{equation}
	\left(\der_X\omega\right)(Y) := X\left(\omega(Y)\right) - \annihprod{\lderb X Y,\omega},
\end{equation}
where $\discformsk 1 M$ is the set of sections of $T^*M$ smooth at the points of $M$ where the signature is constant.
\end{definition}

\begin{definition}
\label{def_cov_der_smooth}
Let $(M,g)$ be a radical-stationary semi-Riemannian manifold. We define
the following vector spaces of differential forms having smooth covariant derivatives:
\begin{equation}
	\srformsk 1 M = \{\omega\in\annihforms M|(\forall X\in\fivect M)\ \der_X\omega\in\annihforms M\},
\end{equation}
\begin{equation}
	\srformsk k M := \bigwedge^k{}_M\srformsk 1 M.
\end{equation}
\end{definition}

\begin{definition}
\label{def_riemann_curvature}
We define the \textit{Riemann curvature tensor} as
\begin{equation}
	R: \fivect M\times \fivect M\times \fivect M\times \fivect M \to \R,
\end{equation}
\begin{equation}
\label{eq_riemann_curvature}
	R(X,Y,Z,T) := (\dera X {\ldera Y}Z - \dera Y {\ldera X}Z - \ldera {[X,Y]}Z)(T)
\end{equation}
for any vector fields $X,Y,Z,T\in\fivect{M}$.
\end{definition}

%--------------------------------------------------------
\subsection{Semi-regular semi-Riemannian geometry}

There is a special type of metric -- which we call \textit{semi-regular} -- which, even though it is degenerate, it is very well behaved. In particular, it allows the definition of smooth covariant derivatives for an important class of differential forms and covariant tensors. It also allows us to construct in a canonical way a smooth Riemann curvature tensor $R_{abcd}$ (unlike $R^a{}_{bcd}$, it can be constructed).

\begin{definition}
\label{def_semi_regular_semi_riemannian}
A \textit{semi-regular semi-Riemannian manifold} is a singular semi-Riemannian manifold $(M,g)$ which satisfies
\begin{equation}
	\ldera X Y \in\srformsk 1 M
\end{equation}
for any vector fields $X,Y\in\fivect{M}$.
\end{definition}

\begin{proposition}
\label{thm_sr_cocontr_kosz}
Let $(M,g)$ be a radical-stationary semi-Riemannian manifold. Then, the manifold $(M,g)$ is semi-regular if and only if for any $X,Y,Z,T\in\fivect M$
\begin{equation}
	\kosz(X,Y,\cocontr)\kosz(Z,T,\cocontr) \in \fiscal M.
\end{equation}
\end{proposition}

\begin{theorem}
\label{thm_riemann_curvature_semi_regular}
Let $(M,g)$ be a semi-regular semi-Riemannian manifold. The Riemann curvature is a smooth tensor field $R\in\tensors 0 4 M$.
\end{theorem}

\begin{example}
\label{s_semi_reg_semi_riem_man_example_diagonal}
An important example of semi-regular metric is provided by taking the metric to be diagonal in a coordinate chart \cite{Sto11a}. The Koszul form, when expressed in coordinates, reduces in fact to Christoffel's symbols of the first kind, which, when the metric is diagonal, are of the form $\pm\frac 1 2\partial_a g_{bb}$. If $g=\sum_a\varepsilon_a\alpha_a^2\de x^a\otimes \de x^a$, where $\varepsilon_a\in\{-1,1\}$, for the metric to be semi-regular is needed that for any $a,b\in\{1,\ldots,n\}$ and $c\in\{a,b\}$ there is a smooth function $f_{abc}\in\fiscal{M}$ so that $\supp{f_{abc}}\subseteq\supp{\alpha_c}$ and
\begin{equation}
\label{eq_diagonal_metric:semireg}
\partial_a\alpha_b^2=f_{abc}\alpha_c.
\end{equation}
Please note that if $c=b$, from $\partial_a\alpha_b^2=2\alpha_b\partial_a\alpha_b$ follows that the function is $f_{abb}=2\partial_a\alpha_b$. This has to satisfy the additional condition that $\partial_a\alpha_b=0$ whenever $\alpha_b=0$. The condition $\supp{f_{abc}}\subseteq\supp{\alpha_c}$ is required because for a manifold to semi-regular it first has to be radical-stationary.
\end{example}

More examples of semi-regular semi-Riemannian manifolds are given in \cite{Sto11b}, where we study the warped product of such manifolds, and conditions which ensure its semi-regularity.

%--------------------------------------------------------
\section{Einstein's equation on semi-regular spacetimes}
The Einstein tensor is usually defined on a semi-Riemannian manifold by the the Ricci tensor and the scalar curvature:
\begin{equation}
\label{eq_einstein_tensor}
	G:=\ric-\frac 1 2 s g
\end{equation}
If the metric is degenerate and radical-stationary, the Ricci tensor and the scalar curvature can be defined as in \cite{Sto11a}. They are smooth as long as the metric doesn't change its signature, but they can become infinite at the points where the metric changes its signature. In \cite{Sto11a} we showed that, if the metric is semi-regular, we can remove the singularity by using instead of the tensorial equation a densitized version.

\begin{definition}
\label{def_semi_reg_spacetime}
A four-dimensional semi-regular semi-Riemannian manifold having the signature $(0,3,1)$ at the points where it is non-degenerate is named \textit{semi-regular spacetime}.
\end{definition}

\begin{theorem}
\label{thm_densitized_einstein}
The Einstein density tensor $G\det g$ on a semi-regular spacetime $(M,g)$ is smooth.
\end{theorem}
If the metric is not degenerate at $p$, we can express the Einstein tensor \eqref{eq_einstein_tensor} using the Hodge $\ast$ operator:
\begin{equation}
\label{eq_einstein_tensor_hodge}
	G_{ab} = g^{st}(\ast R\ast)_{asbt},
\end{equation}
where $(\ast R\ast)_{abcd}$ is the double Hodge dual of $R_{abcd}$ with respect to both the first and the second pairs of indices \cf \eg \cite{PeR87}, p. 234. This can be written explicitly, in terms of the components $\varepsilon_{abcd}$ of the volume form associated to the metric, as
\begin{equation}
	(\ast R\ast)_{abcd} = \varepsilon_{ab}{}^{st}\varepsilon_{cd}{}^{pq}R_{stpq}.
\end{equation}

The volume form can be expressed in coordinates, in terms of the Levi-Civita symbol, by
\begin{equation}
\varepsilon_{abcd} = \epsilon_{abcd}\sqrt{-\det g},
\end{equation}
which allows the rewriting of the Einstein tensor as
\begin{equation}
\label{eq_einstein_tensor_hodge_lc}
	G^{ab} = \frac{g_{kl}\epsilon^{akst}\epsilon^{blpq} R_{stpq}}{\det g},
\end{equation}

At the points where the signature changes, the Einstein tensor is not guaranteed to be smooth and it may become infinite. The tensor density $G^{ab}\det g$ in turn, remains smooth:
\begin{equation}
\label{eq_einstein_tensor_density}
	G^{ab}\det g = g_{kl}\epsilon^{akst}\epsilon^{blpq} R_{stpq}
\end{equation}
in a semi-regular spacetime. This is because $G^{ab}\det g$ is constructed from the smooth Riemann curvature tensor (\cf Theorem \ref{thm_riemann_curvature_semi_regular}), and from the Levi-Civita symbol, which is also smooth. The tensor density $G_{ab}\det g$ obtained by lowering the indices, is smooth too.

These observations allow us to write a \textit{densitized Einstein equation}:
\begin{equation}
\label{eq_einstein:densitized}
	G\det g + \Lambda g\det g = \kappa T\det g,
\end{equation}
or, locally,
\begin{equation}
\label{eq_einstein_idx:densitized}
	G_{ab}\det g + \Lambda g_{ab}\det g = \kappa T_{ab}\det g,
\end{equation}
where $\mathcal G$ and $c$ are Newton's constant and the speed of light, $\kappa:=\frac{8\pi \mathcal G}{c^4}$, and $T$ is the stress-energy tensor.

Since we could construct a differential geometry and a form of Einstein's equation which works with a class of singularities, we may reconsider the thought that general relativity breaks down with necessity at singularities.

%--------------------------------------------------------
\section{Physical laws in general relativity}
The physical laws are described by \textit{evolution equations}.

Einstein's equation \eqref{eq_einstein_idx} shows the relation between the \textit{matter fields} (present through the stress-energy tensor $T_{ab}$), and the \textit{geometry} (represented by Einstein's tensor $G_{ab}:=R_{ab}-\frac 1 2 g_{ab}R$).

To have a well defined evolution, it is normally required that

\begin{enumerate}
	\item
	the quantities involved in Einstein's equation are \textit{non-singular}
	\item
	the spacetime admits a \textit{space-like foliatio}n by a family of Cauchy hypersurfaces.
\end{enumerate}

%--------------------------------------------------------
\subsection{The problem with singularities}

Singularities violate both requirements of a well defined evolution:

\begin{enumerate}
	\item
They make the quantity involved in Einstein's equation become infinite.
	\item
They seem to block the time evolution, because the Cauchy data is partially lost in the singularities.
\end{enumerate}

Therefore, the information seems to no longer be preserved, leading to the \textit{black hole information paradox} \cite{Haw73,Haw76}.

%--------------------------------------------------------
\subsection{Repairing Einstein's equation}

The problem of singular quantities involved in the evolution equations can be approached by a method used by Einstein and Rosen \cite{ER35} (for which they credited Mayer). They suggested to multiply the Riemann and the Ricci tensors by a power of $\det g$, so that all occurrences of $g^{ab}$ in the expression of the Ricci and scalar curvatures are replaced by the adjugate matrix of $g_{ab}$, $\det g g^{ab}$ (\cite{ER35}, p. 74).

Another possibility is to use other appropriate quantities to multiply with. For example, consider the Kerr-Newman solution of Einstein's equations expressed in Boyer-Lindquist coordinates \cite{BL67,Wal84}, p. 313:

\begin{equation}
\label{eq_kerr_newman}
\begin{array}{lll}
\de s^2 &=& \frac{\Delta - a^2\sin^2\vartheta}{\Sigma}\de t^2 - \frac{2a\sin^2\vartheta(r^2 + a^2 - \Delta)}{\Sigma}\de t \de\varphi \\
&&+ \frac{(r^2 + a^2)^2 - \Delta a^2 \sin^2\vartheta}{\Sigma}\sin^2 \vartheta \de\varphi^2 + \frac{\Sigma}{\Delta}\de r^2 + \Sigma\de\vartheta^2,\\
\end{array}
\end{equation}
where
\begin{equation}
\Sigma:=r^2 + a^2\cos^2\vartheta,
\end{equation}
\begin{equation}
\Delta:=r^2 + a^2 + e^2 - 2Mr,
\end{equation}
and $e$ is the electric charge, $a$ the angular momentum, $M$ the mass. This solution contains as particular cases the Schwarzschild, Reissner-Nordstr\"om and Kerr solutions.

We can multiply equation (\ref{eq_kerr_newman}) with the product $\Sigma\Delta$ and obtain a non-singular expression. The equation \eqref{eq_kerr_newman} modified like this will depend on the coordinate system, but it remains non-singular. Of course, if we would want to solve for the geometry (\eg for the metric), we will still obtain infinities, and we may apply the method used by Einstein and Rosen.

While this method may work in replacing the singular data with smooth values, from geometric viewpoint one would be happier with an invariant theory which allows us to deal properly with such singularities, by relying on smooth fields. Our development of singular semi-Riemannian geometry for metric with variable signature presented in the earlier sections makes some steps in this direction, but for the important singularities in general relativity we may need a more powerful tool. This is because the singularities of interest for relativity are not necessarily of degenerate type (as we can see from equation (\ref{eq_kerr_newman})) \footnote{Recent progresses, taking place after this conference, show that the singularities can be made of degenerate type, at least for Schwarzschild and Reissner-Nordstr\"om solutions (see section \S\ref{s_addendum}).}.

%--------------------------------------------------------
\section{Repairing the topology}
Once the equations are repaired, we need to solve the second problem: finding the correct topology at the singular regions. Here usually the things are ambiguous, because we cannot tell which coordinate systems are singular at the singularities and which are not. This freedom allows us to find a nice topology, which allows us to select globally hyperbolic regions from the maximally extended spacetimes, and to find space-like foliations of these regions.

We will exemplify this method on the standard black hole solutions.

%--------------------------------------------------------
\subsection{Penrose-Carter diagrams}

Penrose invented a method to map the entire spacetime on a finite region of a piece of paper. He employed the spherical symmetry of some of the solutions to reduce the number of dimensions to $2$, $r$ and $t$. Then, he applied a conformal mapping to map for example the Minkowski spacetime to a diamond shaped region (Fig. \ref{diamond}). Very often we represent only half of the diagram, corresponding to $r\geq0$.

\image{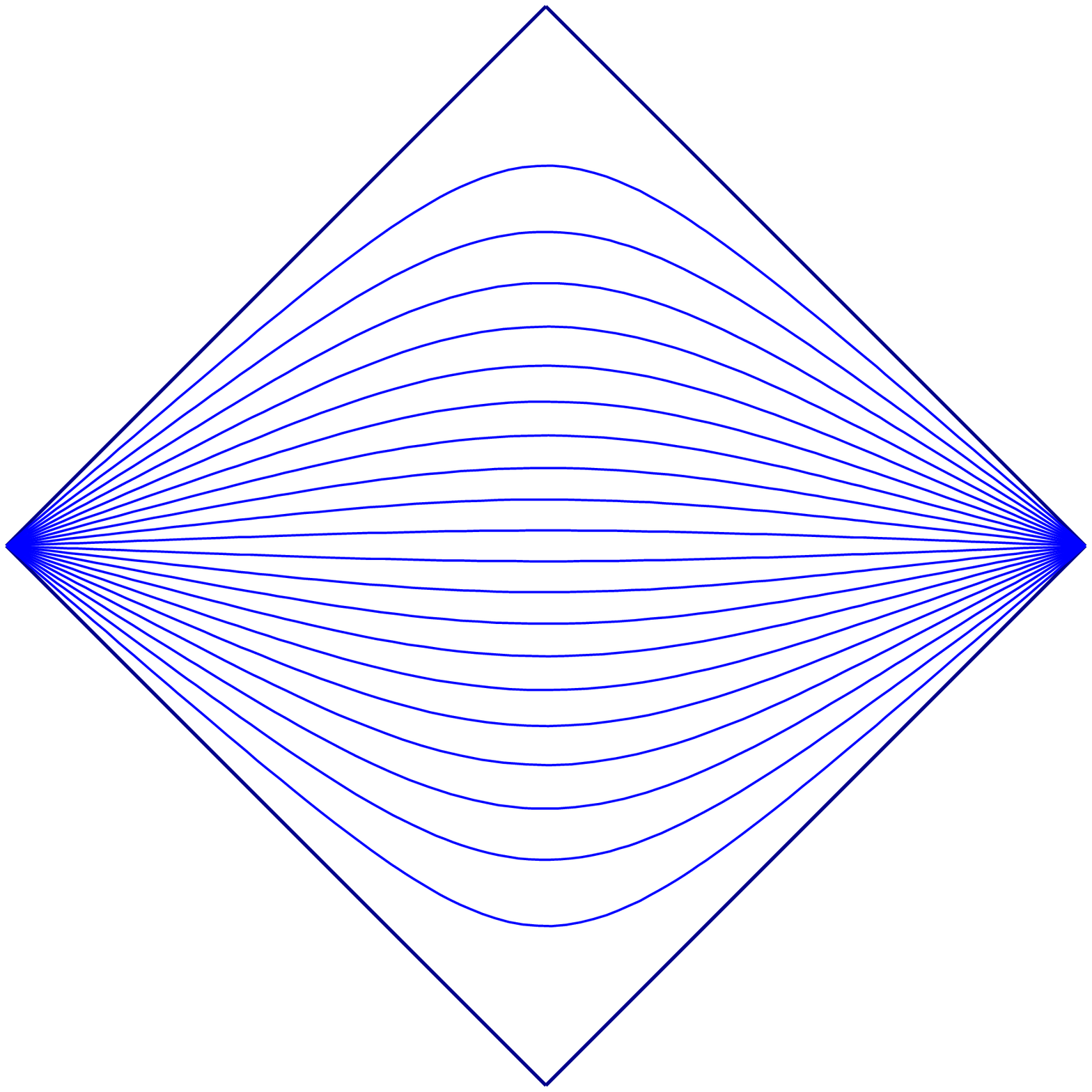}{2.5}{The Penrose-Carter diagram corresponding to the Minkowski spacetime.}

Similar mappings work for other solutions of Einstein's equations.
When Penrose's coordinates are used, they may reveal that the spacetime can be analytically extended. Figure \ref{std-schwarzschild} shows the maximally extended Schwarzschild solution, in Penrose coordinates. It seems to imply that the black hole is paired with a white hole, in the past.

\image{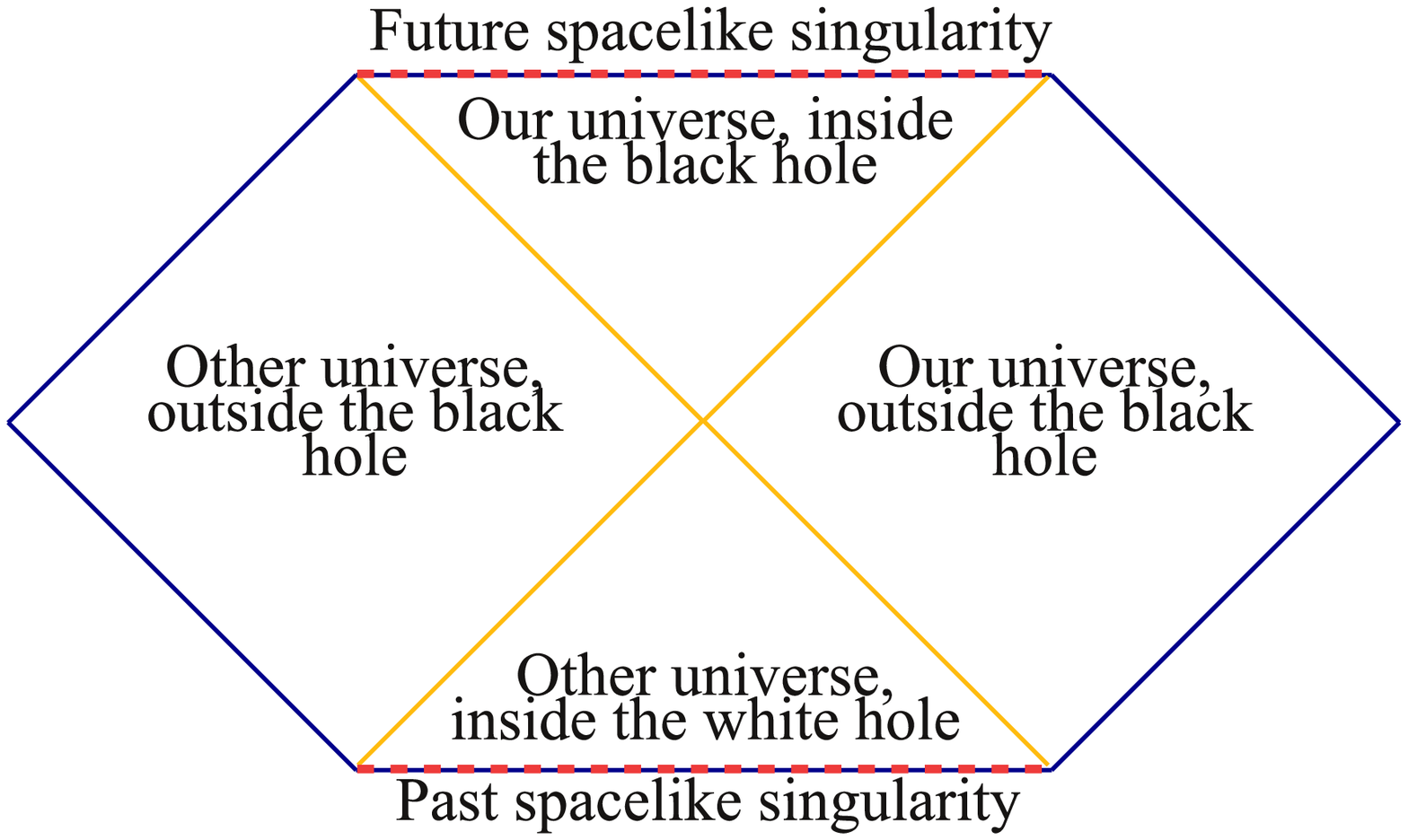}{3}{The maximally extended Schwarzschild solution, in Penrose coordinates.}

Sometimes, the analytical continuation leads to an infinite chain of universes similar to ours. This is the case with the Penrose-Carter diagram of the electrically charged (Reissner-Nordstr\"om) black holes (Fig. \ref{std-rn}). The Penrose-Carter diagrams which describe Kerr (rotating) black holes are very similar to those for the Reissner-Nordstr\"om black holes, except that the symmetry is now axial, and at $r=0$ there is analytical continuation to negative values for $r$.

\image{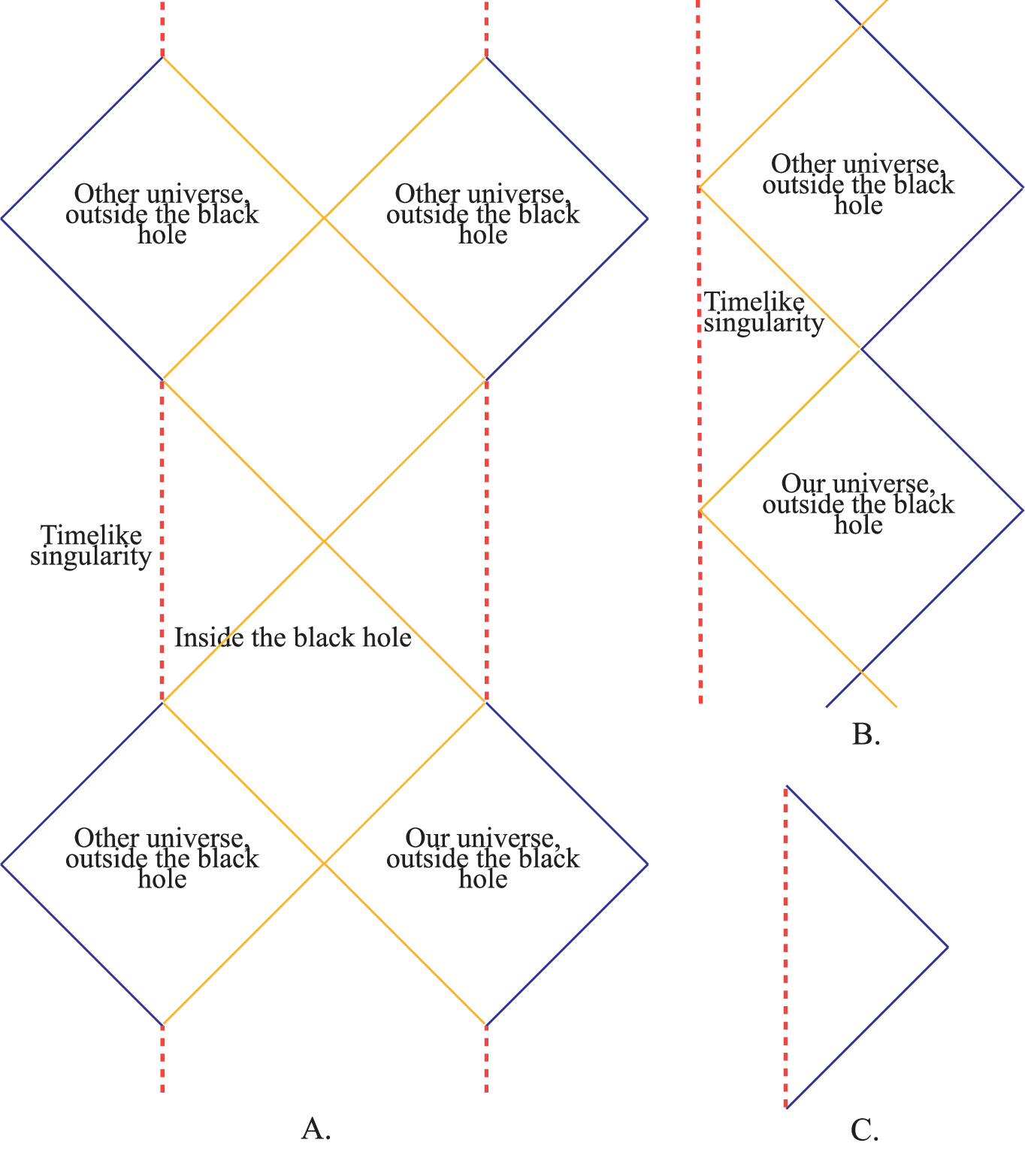}{4}{A. Reissner-Nordstr\"om black holes with $e^2<m^2$. B. Extremal Reissner-Nordstr\"om black holes ($e^2=m^2$). C. Naked Reissner-Nordstr\"om black holes ($e^2>m^2$).}

The maximally extended solution may contain \textit{Cauchy horizons}, so if we want a good foliation of the spacetime, we need to keep only a globally hyperbolic region of it.

%--------------------------------------------------------
\subsection{Schwarz-Christoffel mapping}
The foliations we want can be obtained with the help of the Schwarz-Christoffel mapping. It has a version mapping the strip

\begin{equation}
\label{eq_strip}
\mathcal S:=\{z\in\C|\textnormal{Im}(z)\in[0,1]\}
\end{equation}
to a polygonal region from $\C$, with the help of the formula
\begin{equation}
	\label{eq_s_c_map}
	f(z)=A + C\int^{\mathcal S}\exp\left[\frac\pi 2(\alpha_--\alpha_+)\zeta\right]\prod_{k=1}^n\left[\sinh \frac\pi 2(\zeta-z_k)\right]^{\alpha_k-1}\de\zeta,
\end{equation}
where $z_k\in\partial\mathcal S:=\R\times\{0,i\}$ are the prevertices of the polygon, and $\alpha_-,\alpha_+,\alpha_k$ are the measures of the angles of the polygon, divided by $\pi$ (\cf \eg \cite{dri02}). The ends of the strip, which are at infinite, correspond to the vertices having the angles $\alpha_-$ and $\alpha_+$. The foliation is given by the level curves $\{\textnormal{Im}(z)=\textnormal{const.}\}$.

%--------------------------------------------------------
\subsection{Foliating the maximally extended Schwarzschild solution}

The Schwarzschild black hole solution has the following metric tensor in Schwarzschild coordinates:

\begin{equation}
\label{eq_schw_schw}
\de s^2 = -\(1-\frac{2m}{r}\)\de t^2 + \(1-\frac{2m}{r}\)^{-1}\de r^2 + r^2\de\sigma^2,
\end{equation}
where
\begin{equation}
\label{eq_sphere}
\de\sigma^2 = \de\theta^2 + \sin^2\theta \de \phi^2
\end{equation}
is the metric of the unit sphere $S^2$, $m$ the mass of the body, and the units were chosen so that $c=1$ and $G=1$ (see \eg \cite{HE95} p. {149}).

Figure \ref{hexagon} presents a space-like foliation of the maximally extended Schwarzschild solution.

\image{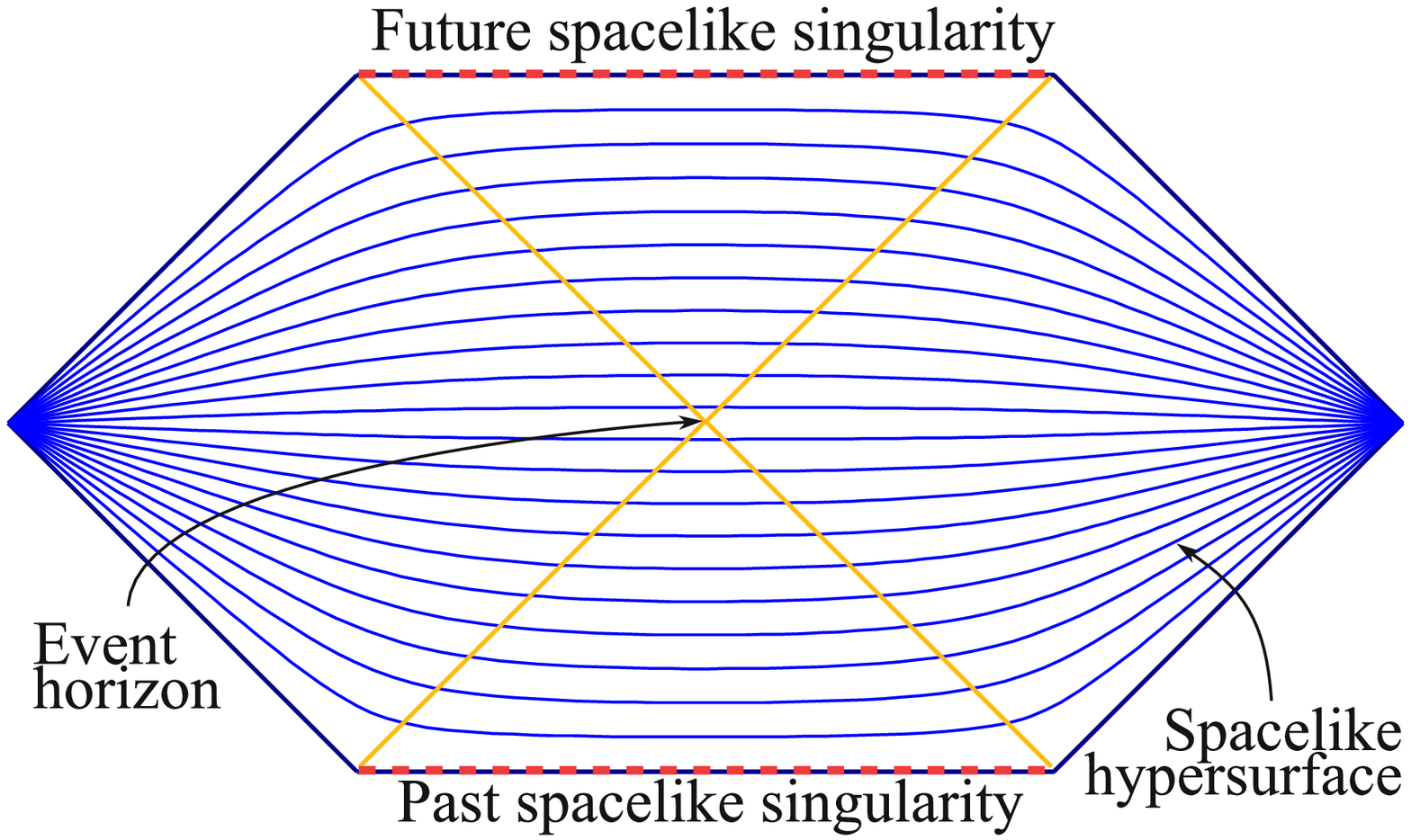}{3}{Space-like foliation of the maximally extended Schwarzschild solution.}

To obtain it, we take the prevertices to be
\begin{equation}
\label{eq_prevertices_hexagon}
	\left(-\infty,-a, a, +\infty, a+i,-a+i\right),
\end{equation}
where $a>0$ is a real number. The angles are respectively
\begin{equation}
\label{eq_angles_hexagon}
\left(\frac {\pi}{2},\frac {3\pi}{4},\frac {3\pi}{4},\frac {\pi}{2},\frac {3\pi}{4},\frac {3\pi}{4}\right).
\end{equation}

We can make a similar foliation, this time without the white hole, if we use the prevertices
\begin{equation}
\label{eq_prevertices_3oo3s}
\left(-\infty,-a, 0, a, +\infty, b+i,-b+i\right),
\end{equation}
where $0<b<a$ are positive real numbers (Fig. \ref{3oo3s}). The angles are respectively
\begin{equation}
\label{eq_angles_3oo3s}
\left(\frac {\pi}{2},\frac {\pi}{2},\frac {3\pi}{2},\frac {\pi}{2},\frac {\pi}{2},\frac {3\pi}{4},\frac {3\pi}{4}\right).
\end{equation}

\image{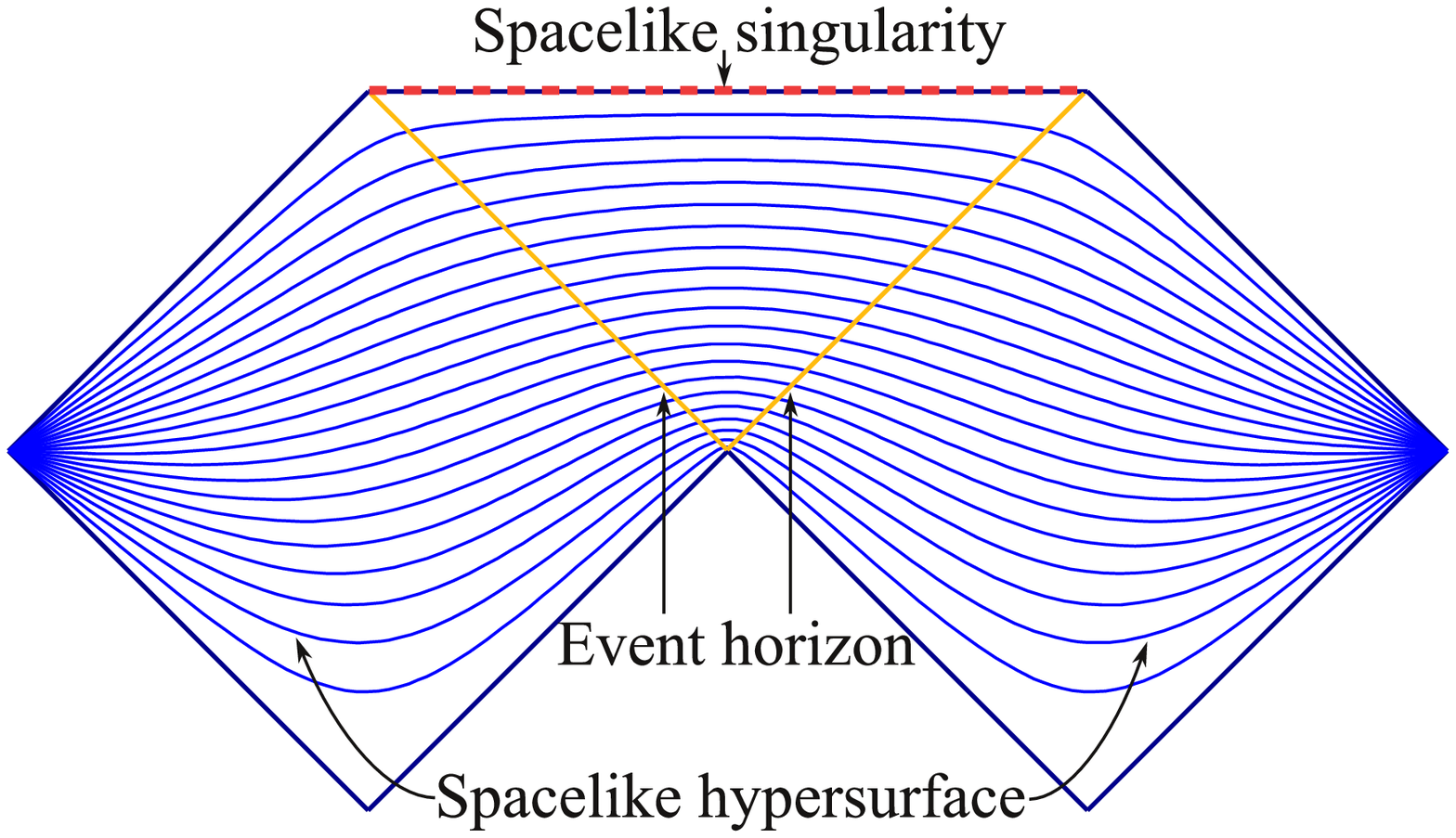}{3}{Space-like foliation of the Schwarzschild solution.}

%--------------------------------------------------------
\subsection{Space-like foliation of the Reissner-Nordstr\"om solution.}

The Reissner-Nordstr\"om solution describes a static, spherically symmetric, electrically charged, non-rotating black hole. It has the following metric:
\begin{equation}
\label{eq_rn_metric}
\de s^2 = -\(1-\frac{2m}{r} + \frac{e^2}{r^2}\)\de t^2 + \(1-\frac{2m}{r} + \frac{e^2}{r^2}\)^{-1}\de r^2 + r^2\de\sigma^2,
\end{equation}
where $e$ is the electric charge of the body, $m$ the mass of the body, and the units were chosen so that $c=1$ and $G=1$ (see \eg \cite{HE95}, p. {156}).

The naked Reissner-Nordstr\"om solution admits a simple foliation, which coincides with that of the Minkowski spacetime (Fig. \ref{diamond-rn-naked}).

\image{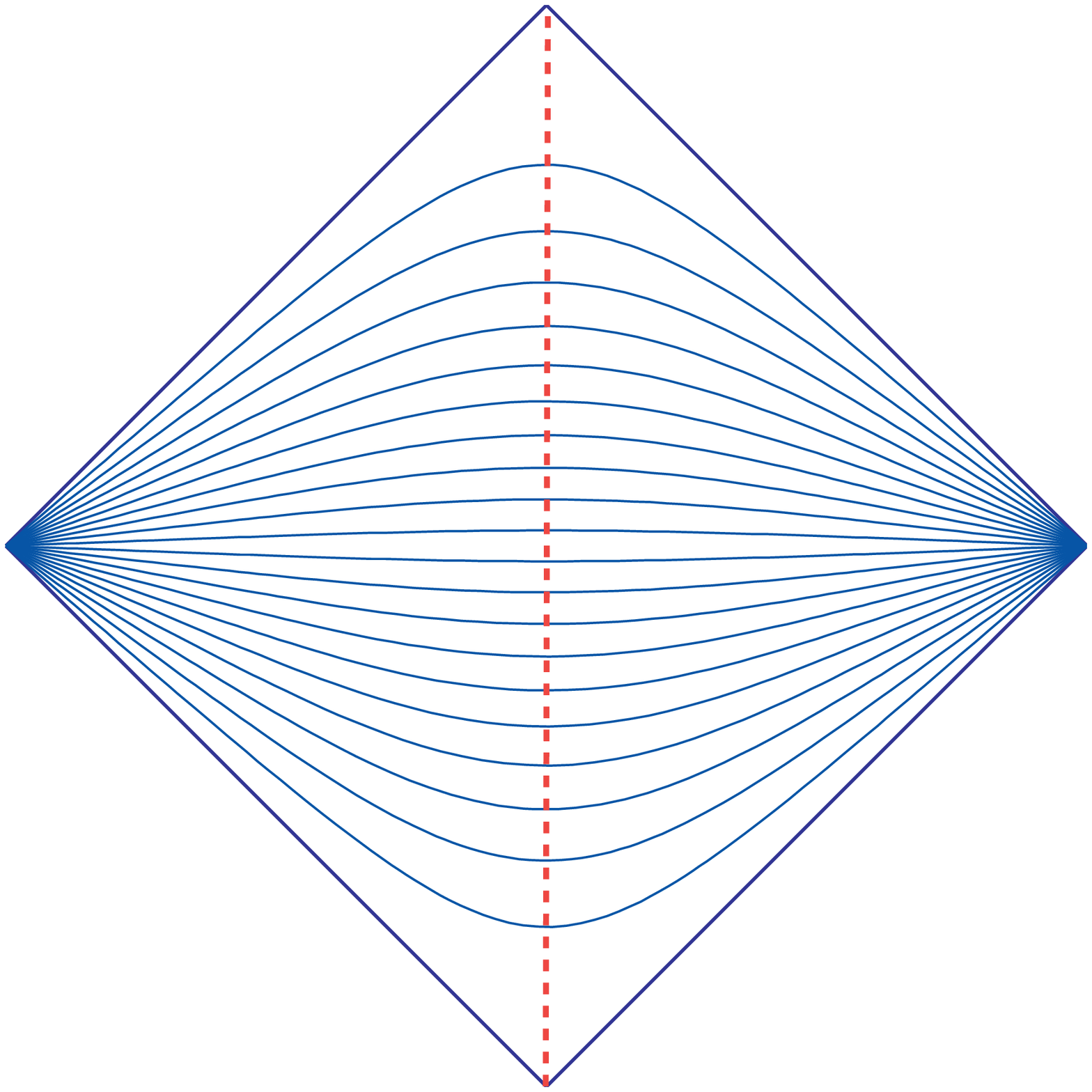}{2.5}{Space-like foliation of the naked Reissner-Nordstr\"om solution ($e^2>m^2$).}

In this case, the prevertices are
\begin{equation}
\label{eq_prevertices_diamond}
\left(-\infty,0, +\infty, i\right),
\end{equation}
and the angles are
\begin{equation}
\label{eq_angles_diamond}
\left(\frac {\pi}{2},\frac {\pi}{2},\frac {\pi}{2},\frac {\pi}{2}\right).
\end{equation}

For the other cases, we take as prevertices of the Schwarz-Christoffel mapping \eqref{eq_s_c_map} the set
\begin{equation}
\label{eq_prevertices_rn-kerr}
\left(-\infty,-a, 0, a, +\infty, i\right),
\end{equation}
where $0<a$ is a positive real number. The angles are respectively
\begin{equation}
\label{eq_angles_rn-kerr}
\left(\frac {\pi}{2},\frac {\pi}{2},\frac {3\pi}{2},\frac {\pi}{2},\frac {\pi}{2},\frac {\pi}{2}\right).
\end{equation}

By a properly chosen $a$, we can obtain a space-like foliation of the non-extremal Reissner-Nordstr\"om solution \ref{up-big}, or of the extremal Reissner-Nordstr\"om solution \ref{up-small-f}.

\image{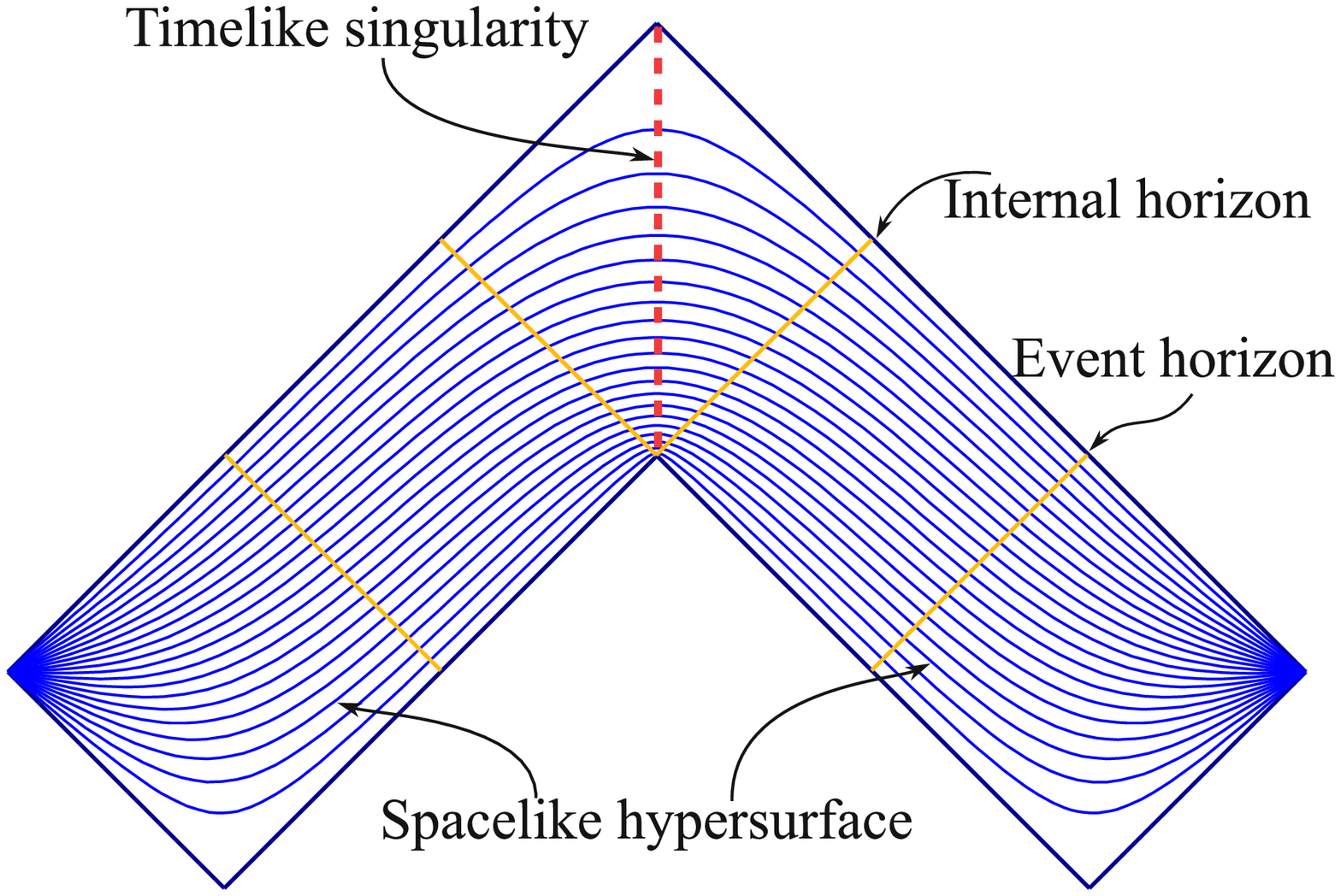}{3.5}{Space-like foliation of the non-extremal Reissner-Nordstr\"om solution ($e^2<m^2$).}

\image{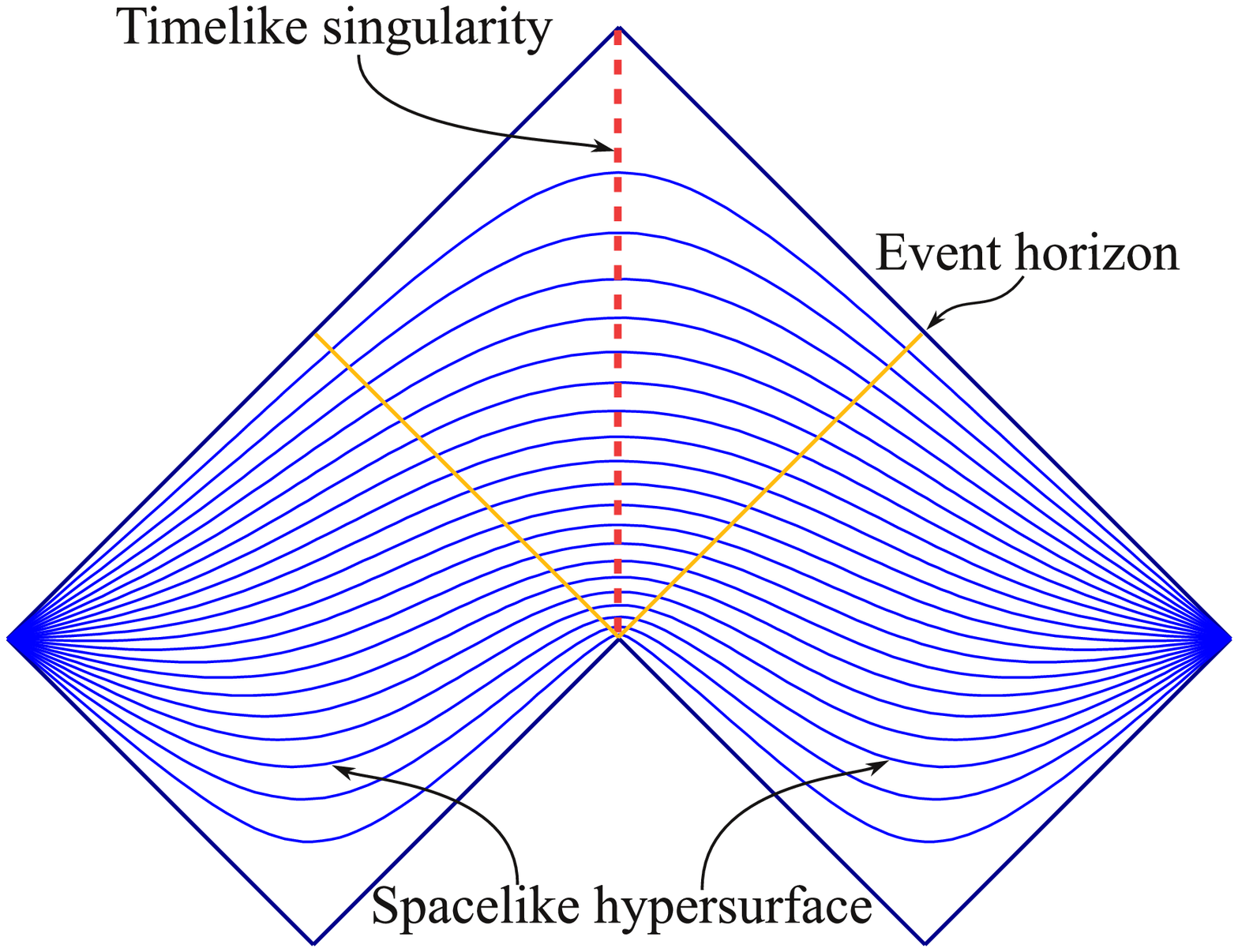}{2.5}{Space-like foliation of the extremal Reissner-Nordstr\"om solution with $e^2=m^2$.}

The case of rotating black holes has many common features with the charged black holes, and the Penrose-Carter diagrams are very similar. But the analysis is more complicated, because it's symmetry is not spherical, only axial. The singularity is a ring, and the analytic extension through the ring leads to closed time-like curves (time machines). We have this case in progress.

%--------------------------------------------------------
\subsection{Non-eternal black holes}

The static solutions are idealized and may not actually represent physically real black holes. In reality, the black holes appear, grow, and possibly evaporate and disappear. There may be primordial black holes, which exist since the beginning of the universe, and black holes which will last until the end of the universe, but they change anyway. In order to have the precise equations of such black holes, we would need to account for all possible types of matter fields which may exist, and the dynamics would be impossible to be solved anyway. But since we shown that the static solutions can be ``cropped'' to spacetimes with trivial topology, we can glue together, in principle, such solutions and obtain at least qualitative descriptions of spacetimes containing black holes with finite life span.

From this viewpoint, we can add before and/or after a spacetime slice as those described above, a spacetime slice without singularities, and model more general black holes. This method relies on the wide flexibility offered by the fact that the Penrose-Carter diagrams are conformal, and one such diagram represents in fact an infinity of possible metrics, which are conformally equivalent to the diagram. We will base this procedure on generalizing the Schwarzschild solution, and the Reissner-Nordstr\"om solution, as prototypes for black holes with space-like, respectively time-like singularities.

%--------------------------------------------------------
\subsubsection{Non-primordial, eternal black holes with space-like singularity}

A space-like foliation of a non-primordial, black hole which continues to exist forever, having a space-like singularity, is shown in Figure \ref{superman}.

\image{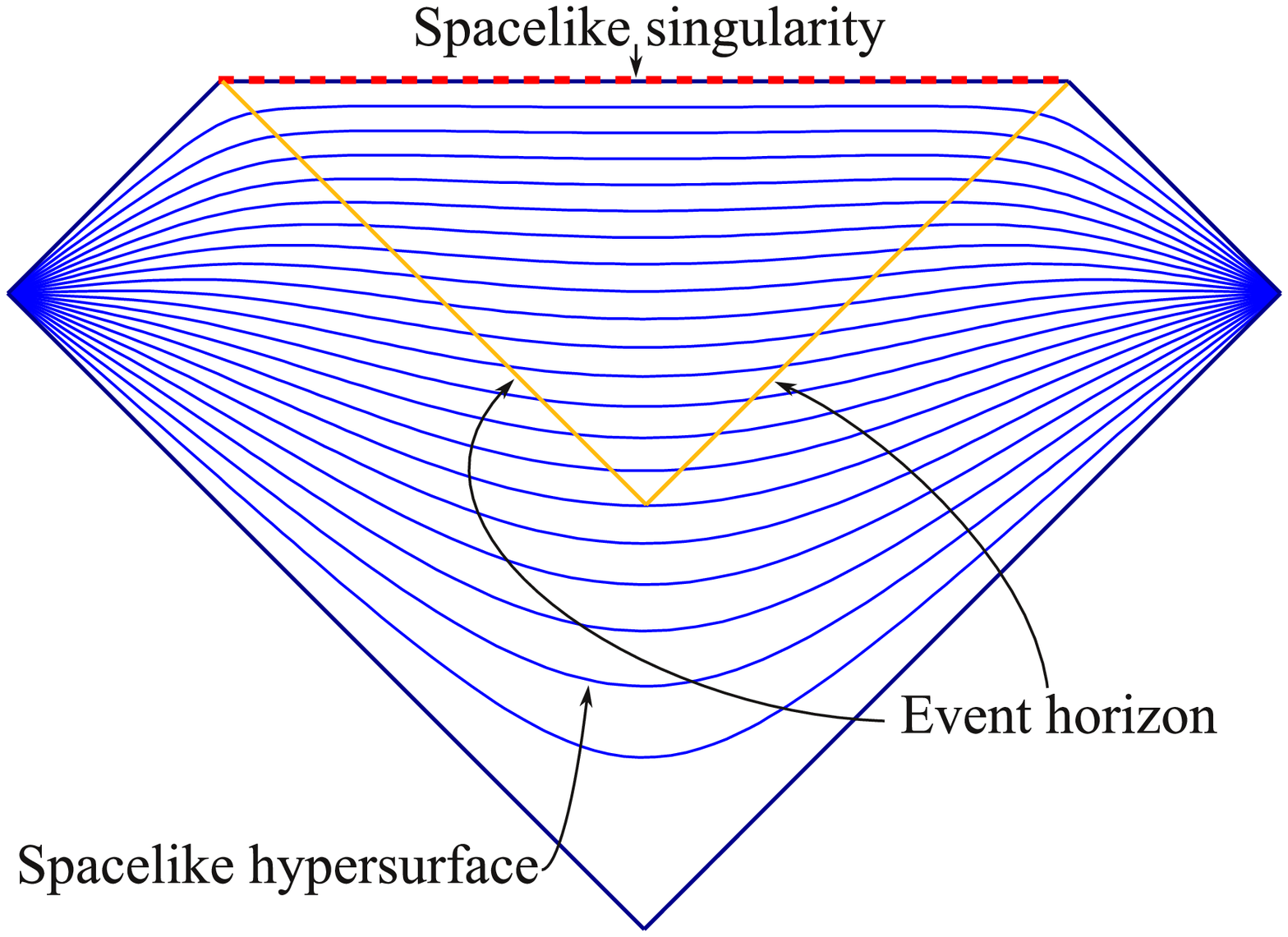}{3}{The space-like foliation for a non-primordial, eternal black hole with space-like singularity.}

The prevertices are given by the set
\begin{equation}
\label{eq_prevertices_superman}
\left(-\infty,0, +\infty, a+i, -a+i\right),
\end{equation}
where $0<a$ is a positive real number. The angles are respectively
\begin{equation}
\label{eq_angles_superman}
\left(\frac {\pi}{2},\frac {\pi}{2},\frac {\pi}{2},\frac {3\pi}{4},\frac {3\pi}{4}\right).
\end{equation}

%--------------------------------------------------------
\subsubsection{Primordial, non-eternal black holes with space-like singularity}

The space-like foliation for primordial black hole which evaporates after a finite time is represented in Fig. \ref{up-small-s}. The prevertices and the angles are the same as in the Reissner-Nordstr\"om solution, as in equations \eqref{eq_prevertices_rn-kerr} and \eqref{eq_angles_rn-kerr}.

\image{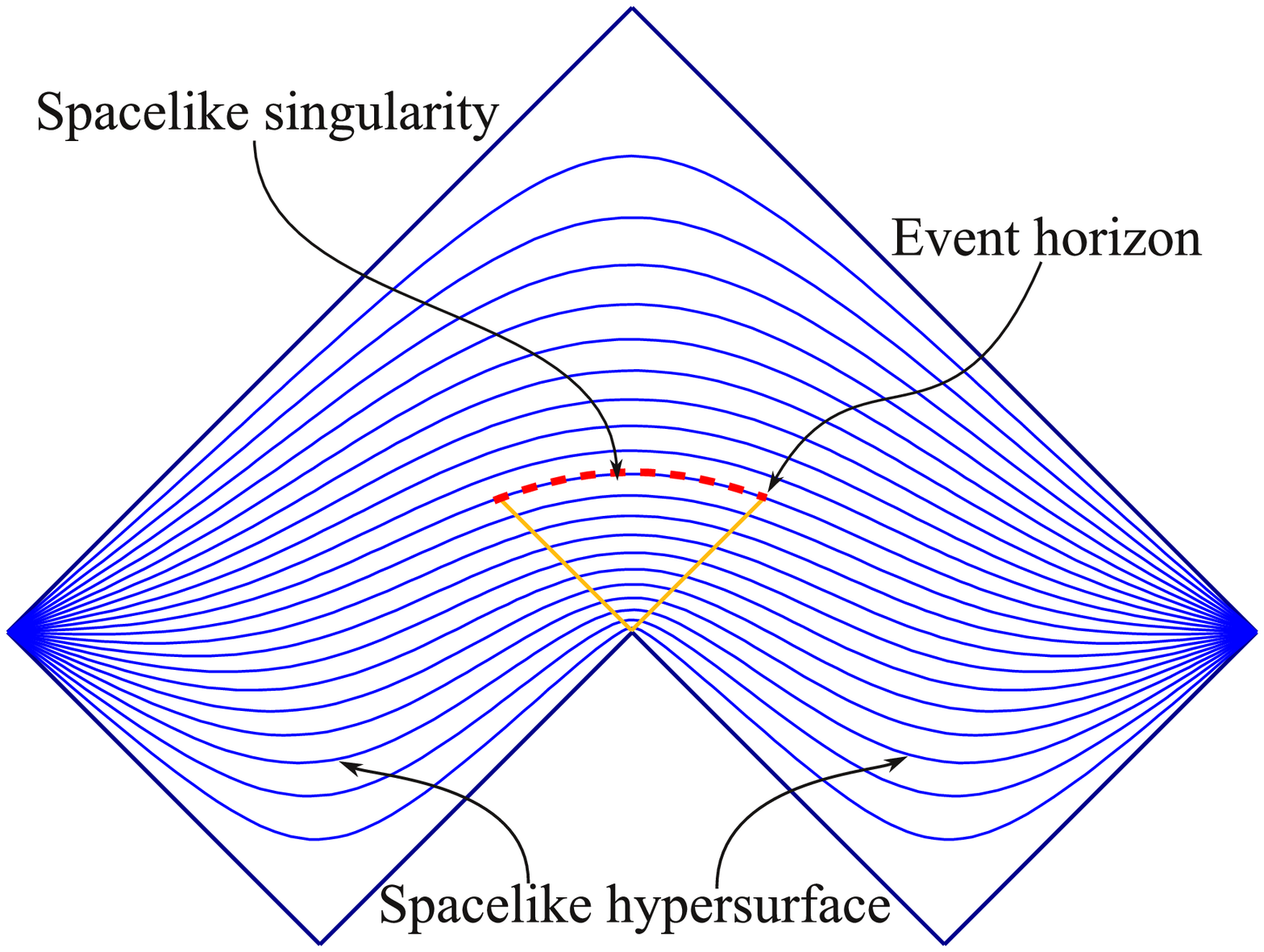}{2.5}{The space-like foliation for a primordial, non-eternal black hole with space-like singularity.}

%--------------------------------------------------------
\subsubsection{Non-primordial, non-eternal black holes with space-like singularity}

In this case (Fig. \ref{diamond-s}), the prevertices are like those in the equation \eqref{eq_prevertices_diamond}, and the angles are the same as in \eqref{eq_angles_diamond}.

\image{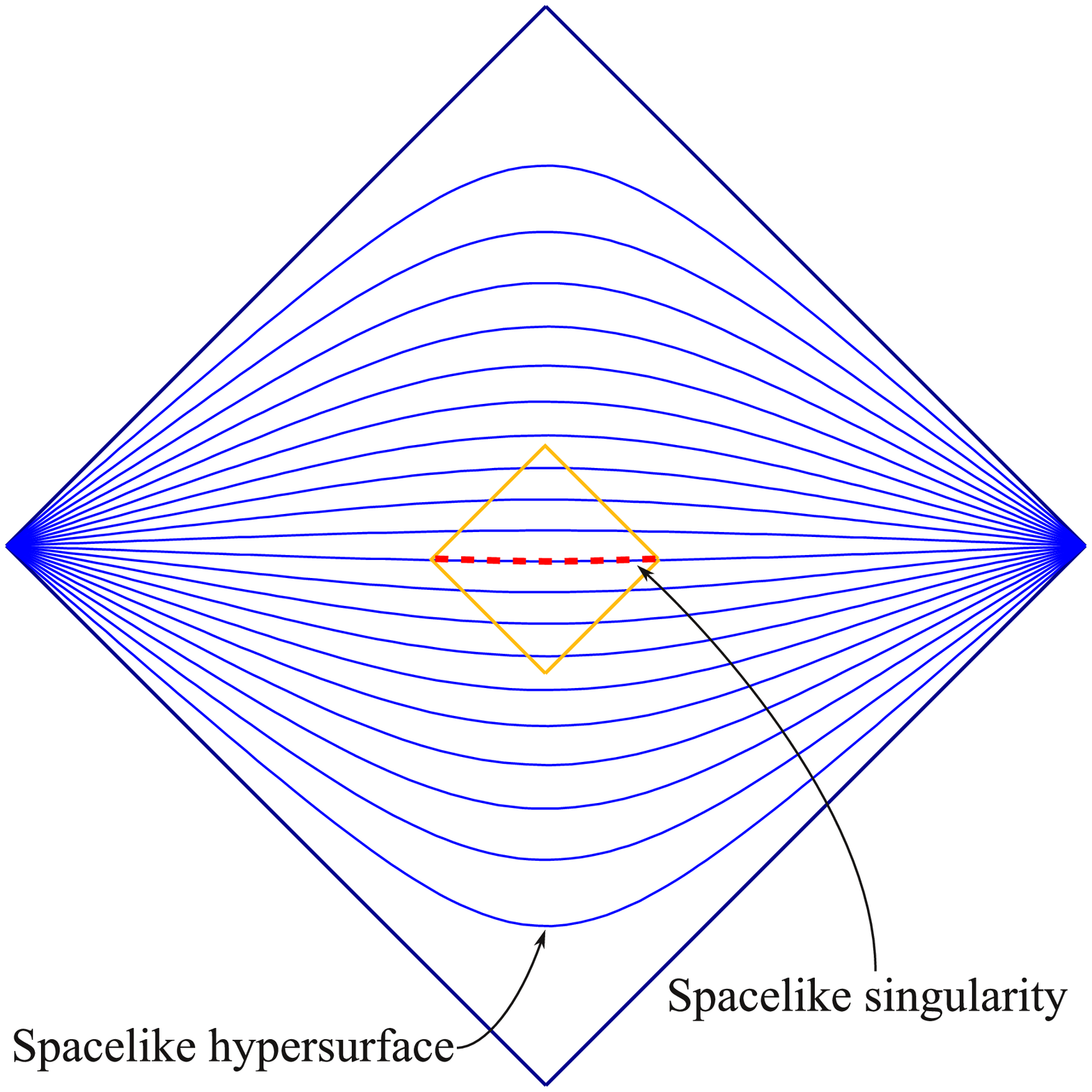}{2.5}{The space-like foliation for a non-primordial, non-eternal black hole with space-like singularity.}

%--------------------------------------------------------
\subsubsection{Primordial evaporating black hole with time-like singularity}

For obtaining the space-like foliation of such a black hole, we take as prevertices and  angles those from the Reissner-Nordstr\"om solution, as in equations \eqref{eq_prevertices_rn-kerr} and \eqref{eq_angles_rn-kerr} (see Fig. \ref{up-small}).

\image{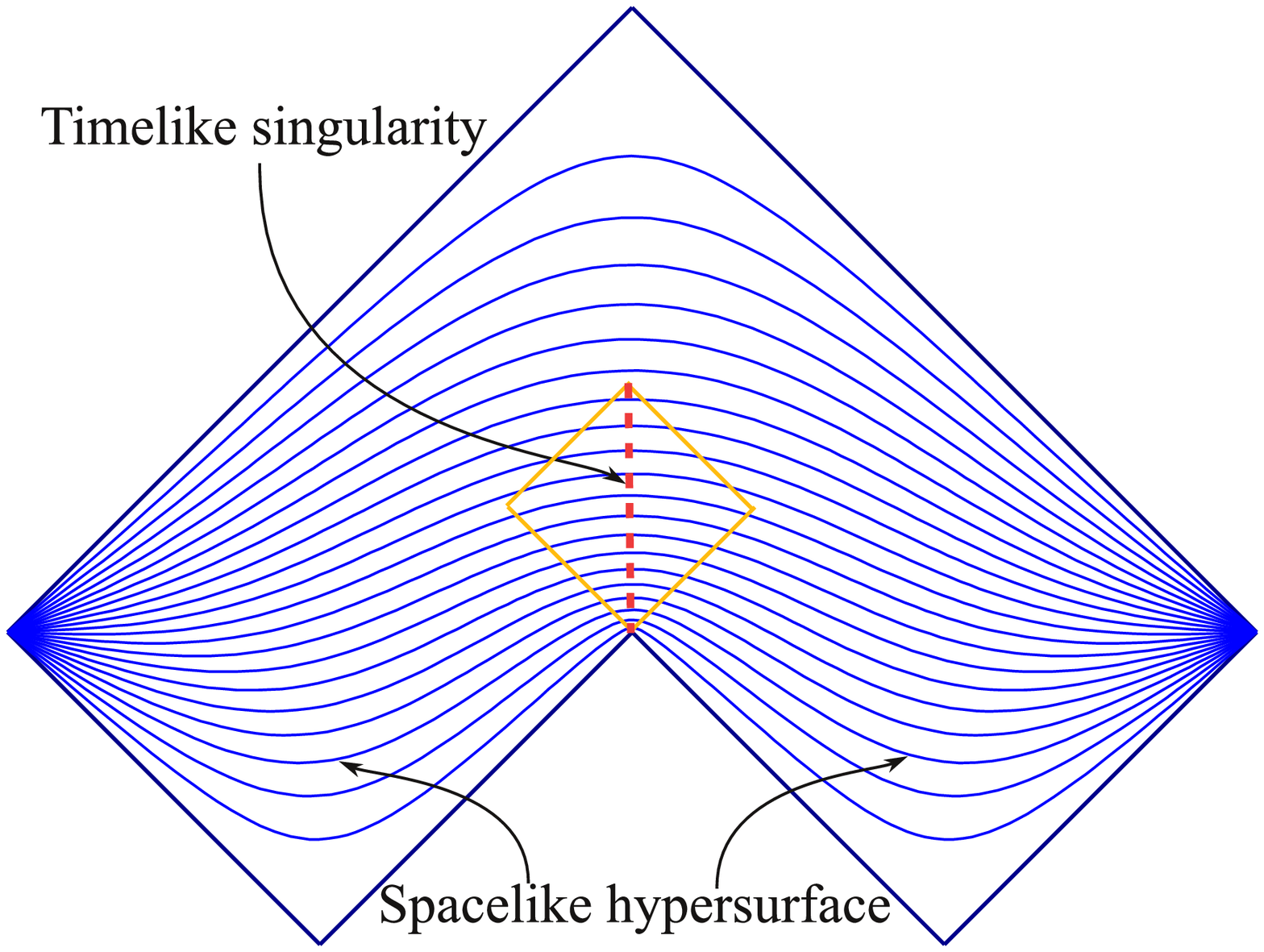}{2.5}{Primordial evaporating black hole with time-like singularity.}

%--------------------------------------------------------
\subsubsection{Non-primordial non-evaporating black hole with time-like singularity}

The prevertices are those from the equation \eqref{eq_prevertices_diamond}, and the angles are the same as in \eqref{eq_angles_diamond}.
In this case, the foliation is the same as for the Minkowski spacetime (Fig. \ref{diamond-t-f}).

\image{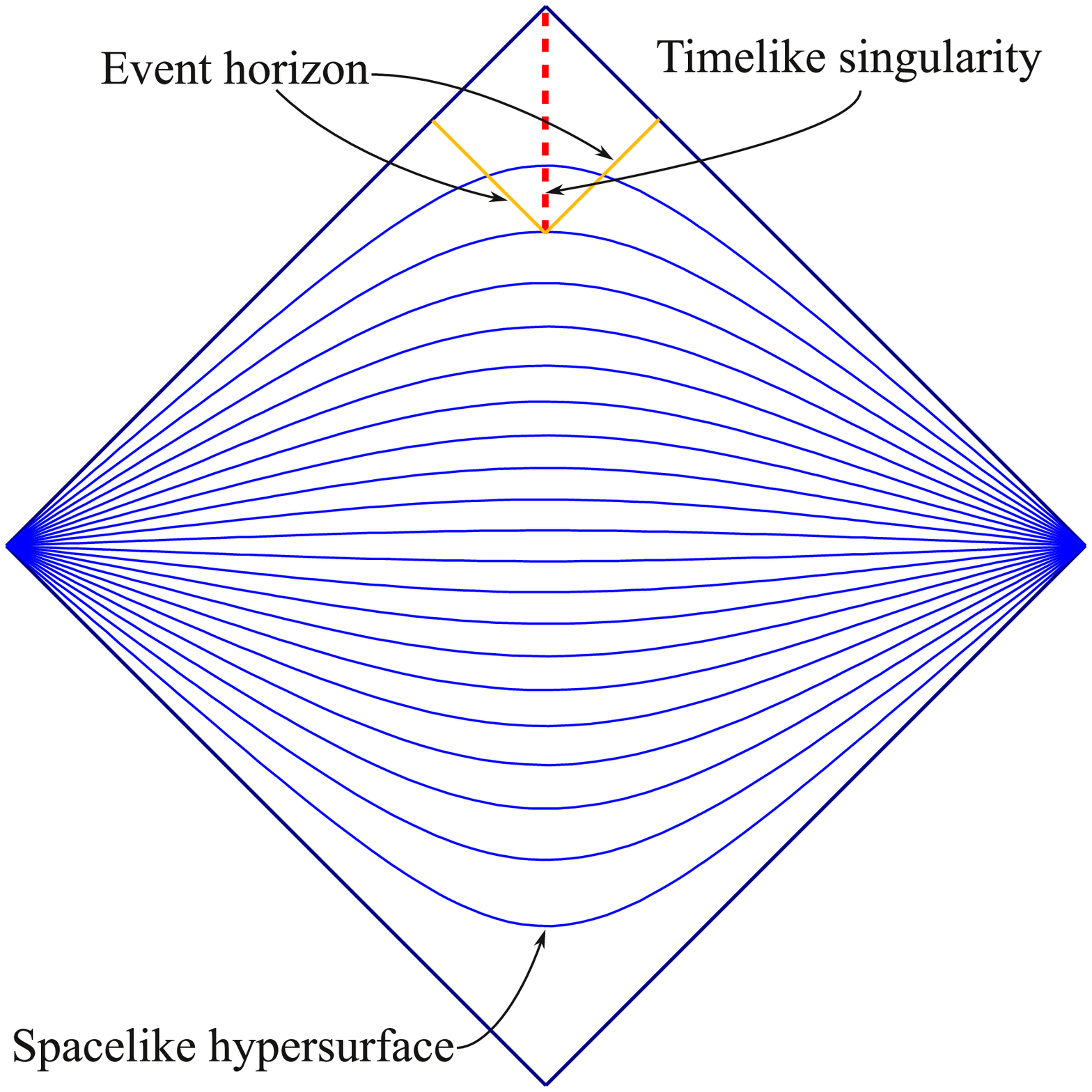}{2.5}{Non-primordial non-evaporating black hole with time-like singularity.}

%--------------------------------------------------------
\subsubsection{Non-primordial evaporating black hole with time-like singularity}

The prevertices are again identical to those from the equation \eqref{eq_prevertices_diamond}, and the angles are the same as in \eqref{eq_angles_diamond} (Fig. \ref{diamond-t}).

\image{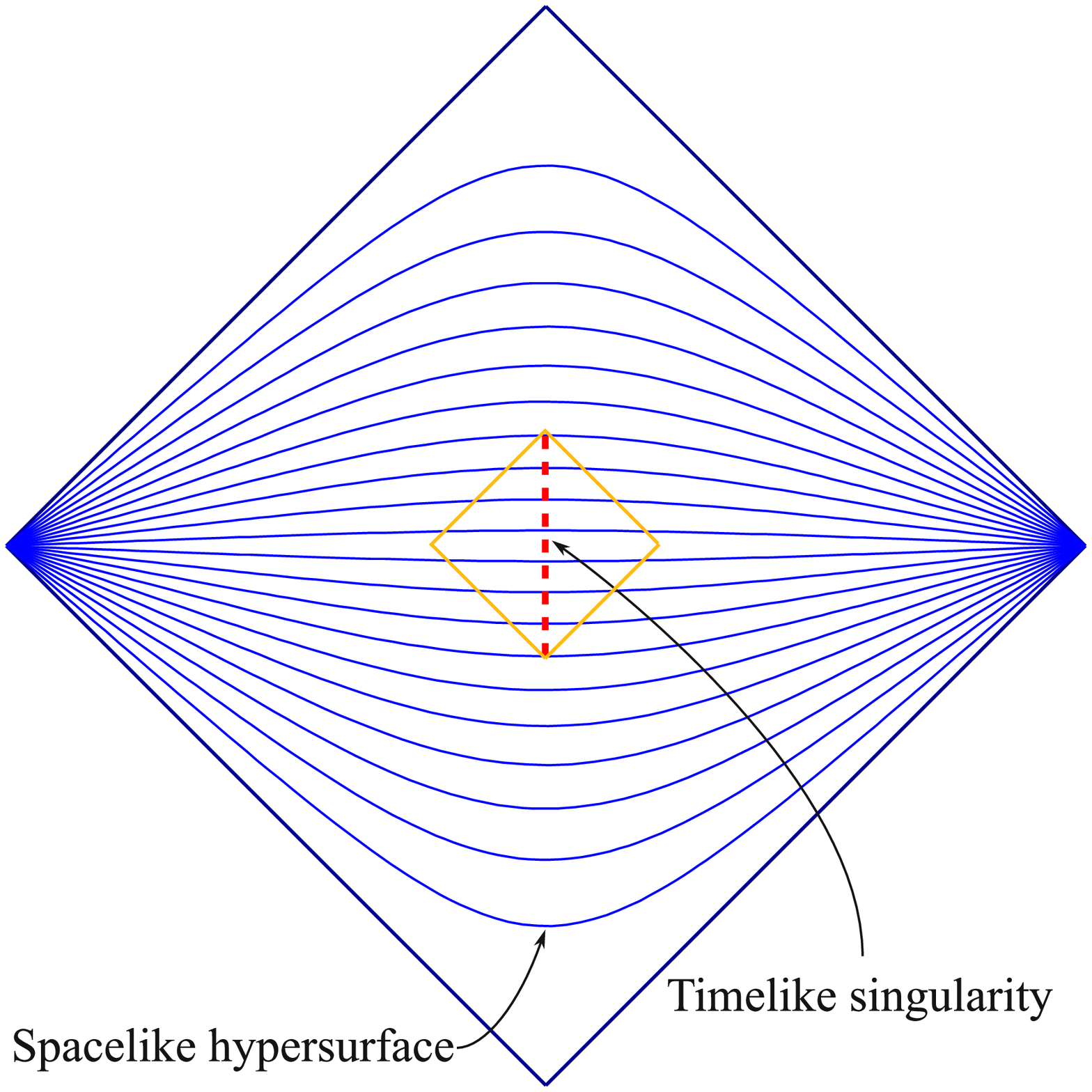}{2.5}{Non-primordial evaporating black hole with time-like singularity.}

%--------------------------------------------------------
\section{Conclusion}

This presentation shows that in important cases we can
\begin{enumerate}
	\item 
make a reformulation of Einstein's equations, in order to express them in terms of finite quantities
	\item 
find an appropriate \textit{foliation} of the spacetime into space-like hypersurfaces
	\item 
choose an appropriate \textit{extension of the spacetime at singularities}, so that the topology of the space-like hypersurfaces of the foliations is preserved,
\end{enumerate}
restoring thus the time evolution at singularities. Consequently, the Cauchy data is preserved, and the information loss is avoided. This shows that the construction of quantum field theories in curved spacetimes with singularities is not necessarily forbidden (\cite{HP96}, p. 9), and the unitarity is not necessarily violated.

%--------------------------------------------------------
\section{For a singular semi-Riemannian geometry}

Once we allow the metric to become degenerate or singular, the invariants like covariant derivative and curvature become undefined. The quantities which replace them if we multiply the equations with $\det g$ or other factors cannot give them this meaning. Maybe a new kind of geometry is required to restore the ideas of covariant derivative and curvature for singular metrics. We did a modest step in this direction, by developing a singular semi-Riemannian geometry for a class of metrics which can become degenerate \cite{Sto11a,Sto11b}, but the more general case of singular metrics remains to be researched.

%--------------------------------------------------------
\section{Short Addendum}
\label{s_addendum}

Actually, the singularities of the standard black hole solutions can be made of degenerate type. Since this conference took place, we made important progress in the direction of dealing with the singularities in a geometric and invariant way. We showed that we can change the coordinates in a way which allows us to extend analytically the Schwarzschild and the Reissner-Nordstr\"om solutions at the singularity and beyond. The Schwarzschild solution turned out to be semi-regularizable \cite{Sto11e}, and the Reissner-Nordstr\"om solution could be analytically extended so that its singularity became of degenerate type only \cite{Sto11f}. This can be viewed as resemblant to Eddington-Finkelstein coordinates, which remove the apparent singularity at the event horizon. In both cases, a singularity is removed by a singular coordinate change, whose singularity overlaps with the singularity of the metric. In our case, the obtained metric is degenerate. One may worry about the fact that this depends on the particular coordinates, but in fact the metric's property of being singular or not is invariant only under non-singular coordinate changes.
We have similar results for the Kerr metric in progress.

%\bibliographystyle{amsalpha}
%\bibliography{../../../bib/sing-gr_bib}
\providecommand{\bysame}{\leavevmode\hbox to3em{\hrulefill}\thinspace}
\providecommand{\MR}{\relax\ifhmode\unskip\space\fi MR }
\providecommand{\MRhref}[2]{%
  \href{http://www.ams.org/mathscinet-getitem?mr=#1}{#2}
}
\providecommand{\href}[2]{#2}

\end{document}